\documentclass[a4paper,leqno,twoside]{article}
\usepackage{latexsym,amssymb,amsmath, amsthm}
\usepackage[affil-it]{authblk}

\newtheorem{theo}{Theorem}[section]
\newtheorem{coro}[theo]{Corollary}
\newtheorem{lemm}[theo]{Lemma}
\newtheorem{prop}[theo]{Proposition}

\newtheorem{rema}[theo]{Remark}

\author{Ivan D. Chipchakov\\
Institute of Mathematics and Informatics, Bulgarian Academy\\ of 
Sciences, Acad. G. Bonchev Str., bl. 8, 1113, Sofia, Bulgaria}
\title{Notes on polynomials $(1 + X)^{n} + (-1) ^{n}(X ^{n} + 1)$ 
concerning the regularity problem for symmetric power sums in 3 
variables}
\date{}
\begin{document}
\maketitle

\begin{abstract}
Let $K$ be a field and $f _{n}(X) = (X + 1) ^{n} + (-1) ^{n}(X ^{n} + 
1) \in K[X]$, for each $n \in \mathbb N$. This note shows that the 
polynomials $f _{m}(X)$ and $f _{m'}(X)$ are relatively prime, for 
some distinct indices $m$ and $m ^{\prime}$ at most equal to $100$, 
if and only if the product $mm ^{\prime }$ is divisible by $6$.  
\end{abstract}


\par
\section{\bf Introduction}
\par
\medskip
Let $f _{n}(X) = (1 + X) ^{n} + (-1)^{n}(X ^{n} + 1)$, for each $n 
\in \mathbb N$. Throught this note, we assume that $f _{n}(X)$, $n 
\in \mathbb N$, are defined over a field $K$ of characteristic zero. 
If the order $n$ of $f _{n}(X)$ is an even number, then the degree 
deg$(f _{n})$ and the leading term of $f _{n}(X)$ are equal to $n$ 
and $2$, respectively; when $n$ is odd, deg$(f _{n})$ and the leading 
term of $f _{n}(X)$ are equal to $n - 1$ and $n$, respectively. In 
addition, it can be easily verified that $f _{n}(X)$ is divisible by the 
polynomial $X(X + 1)$, i.e. $f _{n}(0) = f _{n}(-1) = 0$, if and only 
if $n$ is odd. Similarly, one obtains by straightforward calculations 
that the polynomial $X ^{2} + X + 1$ divides $f _{n}(X)$ if and only 
if $n$ is not divisible by $3$. These observations prove the 
left-to-right implication in the following question:
\par
\medskip
(1) Find whether $f _{m}(X)$ and $f _{n}(X)$ are relatively prime, 
for a pair of distinct positive integers $m$ and $n$, if and only if 
$mn$ is divisible by $6$.
\par
\medskip\noindent
An affirmative answer to (1) would prove the following conjecture in 
the special case where $a = 1$:
\par
\medskip
(2) Let $a$, $b$ and $c$ be a sequence of pairwise distinct positive 
integers with gcd$(a, b, c) = 1$. Then the symmetric polynomials $X 
_{1} ^{a} + X _{2} ^{a} + X _{3} ^{a}$, $X _{1} ^{b} + X _{2} ^{b} + 
X _{3} ^{b}$ and $X _{1} ^{c} + X _{2} ^{c} + X _{3} ^{c}$ form a 
regular sequence in the polynomial ring $K[X _{1}, X _{2}, X _{3}]$ 
in three algebraically independent variables over the field $K$ if 
and only if the product $abc$ is divisible by $6$.
\par
\medskip\noindent
Let us note that a set of $\sigma $ homogeneous polynomials in 
$\sigma $ algebraically independent variables is a regular sequence, 
if the associated polynomial system has only the trivial solution 
$(0, \dots , 0)$. Conjecture (2) has been suggested in \cite{CKW} 
(see also \cite{KM}). 
The purpose of this note is to show that the answer to (1) is 
affirmative, for polynomials of admissible degrees at most equal to 
$100$. Formally, our main result can be stated as follows:
\par
\medskip
\begin{theo}
\label{theo1.1}
Let $m$ and $n$ be distinct positive integers at most equal to $100$, 
and let $K$ be a field with {\rm char}$(K) = 0$. Then the polynomials 
$f _{m}(X), f _{n}(X) \in K[X]$ satisfy the equality {\rm gcd}$\{f 
_{m}(X), f _{n}(X)\} = 1$ if and only if $6 \mid mn$.
\end{theo}
\par
\medskip
It is clearly sufficient to prove Theorem \ref{theo1.1} and to consider 
(1) in the special case where $K$ is the field $\mathbb Q$ of 
rational numbers. Our notation and terminology are standard, and 
missing definitions can be found in \cite{L}.
\par
\medskip
\section{\bf Preliminaries}
\par
\medskip
This Section begins with a brief account of some properties of the 
polynomials $f _{m}(X)$, $f _{n}(X)$, where $m$ and $n$ are distinct 
positive integers. These properties are frequently used in the sequel 
without an explicit reference. Some of the most frequently used facts 
can be presented as follows:
\par
\medskip
(2.1) (a) Any complex root $\alpha _{n}$ of $f _{n}(X)$, for a given 
$n \in \mathbb N$, satisfies the following: $2\alpha _{n}$ is an 
algebraic integer, provided that $n$ is even; $n\alpha _{n}$ is an 
algebraic integer in case $n$ is odd;
\par
(b) If $m$ and $n$ are positive integers of different parity, 
then the common complex roots of $f _{m}(X)$ and $f _{n}(X)$ (if any) 
are algebraic integers;
\par
(c) $0$ and $1$ are simple roots of $f _{n}(X)$, provided that $n \in 
\mathbb N$ and $n$ is odd; the same applies to the reduction of $f 
_{n}(X)$ modulo $2$;
\par
(d) Given an integer $n >0$ and a primitive cubic root of unity 
$\rho _{3}$ (lying in the field $\mathbb C$ of complex numbers), we have 
$f _{n}(\rho _{3}) = 0$ if and only if $n$ is not divisible by $3$.
\par
\medskip
\noindent
Note also that polynomial $f _{n}(X)$ has no real root, for any even 
integer $n$.
\par
\medskip
(2.2) The roots of $f _{p}(X)$ are algebraic integers, for every prime $p$. 
\par
\medskip
Next we include a list of the polynomials $f _{n}(X)$, 
for some small values of $n$:
\par
\medskip
(2.3) (a) $f _{10}(X) = (X ^{2} + X + 1) ^{2}g _{10}(X)$, where $g _{10}(X) = 2X 
^{6} + 6X ^{5} + 27X ^{4} + 44X ^{3} + 27X ^{2} + 6X + 2$;
\par
(b) $f _{9}(X) = 3X(X + 1)g_{9}(X)$, where
\par\noindent
$g _{9}(X) = 3X ^{6} + 9X ^{5} + 
19X ^{4} + 23X ^{3} + 19X ^{2} + 9X + 3$;
\par
(c) $f _{8}(X) = 2(X ^{2} + X + 1)g_{8}(X)$, where
\par\noindent
$g _{8}(X) = X ^{6} + 3X ^{5} + 10X ^{4} + 15X ^{3} + 10X ^{2} + 3X + 1$;
\par
(d) $f _{7}(X) = 7X(X + 1)(X ^{2} + X + 1) ^{2}$;
\par
(e) $f _{6}(X) = 2X ^{6} + 6X ^{5} + 15X ^{4} + 20X ^{3} + 15X ^{2} + 6X + 
2$;
\par
(f) $f _{5}(X) = 5X(X + 1)(X ^{2} + X + 1)$;
\par
(g) $f _{4}(X) = 2(X ^{4} + 2X ^{3} + 3X ^{2} + 2X + 1) = 2(X ^{2} + X + 
1) ^{2}$;
\par
(h) $f _{3}(X) = 3X(X + 1)$;
\par
(i) $f _{2}(X) = 2(X ^{2} + X + 1)$.
\par
\medskip
The following lemma presents some well-known properties of Newton's 
binomial coefficients that are frequently used in the sequel:
\par
\medskip
\begin{lemm}
\label{lemm2.1} 
Assume that  $p$ is a prime number, and $n, s$ are positive integers, 
such that $s$ does not divide $p$. Then the binomial coefficients 
$\binom{p^{n}}{j}$, $j = 1, \dots , p ^{n} - 1$, are divisible by $p$; 
$\binom{p ^{n}s}{p ^{n}}$ is not divisible by $p$.
\end{lemm}
\par
\medskip
\begin{proof} The assertion is obvious if $n = 1$, so we assume 
further that $j \ge p$ (and $n \ge 2$). Suppose first that $j$ is not 
divisible by $p$ and denote by $y$ the greatest integer divisible by 
$p$ and less than $j$. It is clear from the definition of Newton's 
binomial coefficients that the maximal power of $p$ dividing 
$\binom{p^{n}s}{j}$ is greater than the maximal power of $p$ dividing 
$\binom{p^{n}s}{y}$; in particular, this ensures that $p^{2}$ divides 
$\binom{p^{n}s}{j}$. Now fix a positive integer $u < p^{n-1}$ and 
denote by $C[p ^{n}s, pu]$ the product of the multiples of the 
numerator of $\binom{p^{n}s}{pu}$ that are divisible by $p$, divided 
by the product of the multiples divisible by $p$ of the denominator 
of $\binom{p^{n}s}{pu}$. 
It is easily verified that $C[p^{n}s, pu] = \binom{p^{n-1}s}{u}$. This 
allows to complete the proof of Lemma \ref{lemm2.1}, arguing by 
induction on $n$. 
\end{proof}
\par
\medskip
\begin{rema}
\label{rema2.2}
It follows from Lemma \ref{lemm2.1}, applied to $p = 2$, that 
$f _{2^{k}3}(X)$ decomposes over the field $\mathbb Q _{2}$ of 
$2$-adic numbers into a product of three irreducible polynomials of 
degree $2^{k}$ each; one of these polynomials lies in the ring 
$\mathbb Z _{2}[X]$ and is $2$-Eisensteinian over the ring $\mathbb Z 
_{2}$ of $2$-adic integers. In view of our irreducibility criterion, see 
Section 3, this means that if $f _{2^{k}3}(X)$ is reducible over 
$\mathbb Q$, then it decomposes into a product of $3$ $\mathbb 
Q$-irreducible polynomials, say $h _{1}(X)$, $h _{2}(X)$ and $h 
_{3}(X)$ (in fact $h _{j}(X)$, $j = 1, 2, 3$, are irreducible even 
over $\mathbb Q_{2}$). More precisely, the 
action of the symmetric group Sym$_{3}$ on the set of roots of $f 
_{2^{k}3}(X)$ induces bijections $y _{j}$, $j = 2, 3$, from the set 
of roots of $h _{1}(X)$ in $\mathbb Q_{2,{\rm sep}}$ onto the set of 
roots of $h _{j}(X)$, for each index $j$. It is therefore clear that if 
gcd$(f _{2^{k}3}(X), f _{n}(X)) \neq 1$, for some $n \in \mathbb N$, 
then $f _{n}(X)$ and $h _{u}(X)$ have a common root $\xi _{u} \in 
\mathbb Q _{2,{\rm sep}}$, for each $u \in \{1, 2, 3\}$. Thus it 
turns out that if gcd$(f _{2^{k}3}(X), f _{n}(X)) \neq 1$, then $f 
_{2^{k}3}(X) \mid f _{n}(X)$; in particular, $f _{n}(X)$ has a complex 
root that is not an algebraic integer.
\end{rema}
\par
\medskip
Lemma \ref{lemm2.1} ca be supplemented as follows:
\par
\medskip
\begin{lemm}
\label{lemm2.3} 
Let $n \in \mathbb N$ and  $p \in \mathbb P$. Then the binomial 
coefficients $\binom{p^{n}}{j}$ are divisible by $p$, provided that 
$j \in \mathbb Z$ and $1 \le j < p^{n}$; in addition, if $n \ge 
2$, then $\binom{p ^{n}}{j}$ is divisible by $p^{2}$ unless $j= 
p^{n-1}j_{0}$, $1 \le j_{0} \le p - 1$.
\end{lemm}
\par
\medskip
\begin{proof} The former part of our assertion is a special case of 
Lemma \ref{lemm2.1}, so we assume further that $n \ge 2$. Suppose 
first that $j$ is not divisible by $p$ and denote by $y$ the greatest 
integer divisible by $p$ and less than $j$. It is clear from the 
definition of Newton's binomial coefficients that the maximal power 
of $p$ dividing $\binom{p^{n}}{j}$ is greater than the maximal 
power of $p$ dividing $\binom{p^{n}}{y}$; in particular, this ensures that 
$p^{2} \mid \binom{p^{n}}{j}$. Now fix a positive integer $u < 
p^{n-1}$ and define $C[p ^{n}, pu]$ as in the proof of Lemma 
\ref{lemm2.1}. It is easily verified that $C[p^{n}, pu] = 
\binom{p^{n-1}}{u}$. This allows to complete the proof of Lemma 
\ref{lemm2.3}, arguing by induction on $n$.
\end{proof}
\par
\medskip
\section{\bf Polynomials of even orders}
\par
\medskip
This Section begins with a criterion for validity of the equality 
gcd$(f _{m}(X), f _{n}(X)) = 1$ in the special case where $m$ is 
divisible by $6$ and $f _{m}(X)$ is irreducible over $\mathbb Q$.   
\par
\medskip
\begin{prop}
\label{prop3.1}
Let $m$ and $n$ be positive integers with $2 \mid m$ and $6 \mid mn$. 
Put $f _{n}(X) = (1 + X)^{n} + (-1)^{n}(X ^{n} + 1)$, and suppose 
that $m < n$ and the polynomial $f _{m}(X) = (1 + x) ^{m} + X ^{m} + 
1$ is $\mathbb Q$-irreducible or $m = 3.2^{k}$, for some $k \in 
\mathbb N$. Then {\rm gcd}$(f _{m}(X), f _{n}(X)) = 1$ except, 
possibly, under the following conditions:
\par
{\rm (a)} $n \equiv 1 ({\rm mod} \ m - 1)$ and $m \equiv n (mod \ 2^{k+1})$, 
where $k$ is the greatest integer for which $2^{k}$ divides $m$;
\par
{\rm (b)} If $m$ is divisible by $4$, then $n/2 \equiv 1 (mod \ m/2 - 1)$.
\end{prop}
\par
\medskip
\begin{proof}
Suppose that gcd$(f _{m}(X), f _{n}) \neq 1$, for some $n \in \mathbb N$. 
This means that $f _{m}(X)$ and $f _{n}(X)$ have a common root $\rho \in 
\mathbb C$. Note further that the irreducibility of $f _{m}(X)$ over 
$\mathbb Q$ and the assumption that $m$ is even indicate that the 
complex roots of $f _{m}(X)$ are not algebraic integers, so it 
follows from (2.1) (b) that $n \equiv 0 \ ({\rm mod} \ 2)$. Observing 
also that $f _{m}(X) \mid f _{n}(X)$ (in $\mathbb Z[X]$), one 
concludes that $f _{m}(1) \mid f _{n}(1)$ (in $\mathbb Z$), i.e. $2 
^{m} + 2 \mid 2 ^{n} + 2$. It is therefore clear that $2 ^{m-1} + 1 
\mid 2 ^{n-1} + 1$, and since $m - 1$ and $n - 1$ are odd, this 
requires that $m - 1 \mid n - 1$, proving the former part of 
Proposition \ref{prop3.1} (a). The rest of the proof of Proposition 
\ref{prop3.1} relies on Lemma \ref{lemm2.1}, which allows to prove 
that $f _{m}(X)$ and $f _{n}(X)$ have unique divisors $\theta 
_{m}(X)$ and $\theta _{n}(X)$, respectively, over the field $\mathbb 
Q _{2}$, with the following properties: $\theta _{m}$ and $\theta 
_{n}$ are $2$-Eisensteinian polynomials over the ring $\mathbb Z 
_{2}$; the degree of $\theta _{m}(X)$ equals the greatest power of 
$2$ dividing $m$, and the degree of $\theta _{n}$ equals the greatest 
power of $2$ dividing $n$. Observing also that $\theta _{m}(X)$ and 
$\theta _{n}(X)$ can be chosen so that their leading terms be equal 
to $1$, one obtains from the divisibility of $f _{n}(X)$ by $f 
_{m}(X)$ that $\theta _{m}(X) = \theta _{n}(X)$. This result 
completes the proof of Proposition \ref{prop3.1} (a). We turn to the 
proof of Proposition \ref{prop3.1} (b), so we suppose further that $4 
\mid m$. In view of Proposition \ref{prop3.1} (a), this means that $4 
\mid n$, which shows that $f _{m}(\sqrt{-1}) = (2\sqrt{-1}) ^{m/2} + 
2 = 2 ^{m/2} + 2$ and $f _{n}(\sqrt{-1}) = 2 ^{n/2} + 2$. As $f 
_{m}(X) \mid f _{n}(X)$, whence $f _{m}(\sqrt{-1}) \mid f 
_{n}(\sqrt{-1})$ (in the ring $\mathbb Z[\sqrt{-1}]$ of Gaussian 
integers), our calculations lead to the conclusion that $2 ^{(m/2)-1} 
\mid 2 ^{(n/2)-1} + 1$ (in $\mathbb Z$). Taking finally into account 
that $(m/2) - 1$ and $(n/2) - 1$ are odd positive integers, one 
obtains that $(m/2) - 1 \mid (n/2) - 1$. Proposition 
\ref{prop3.1} is proved.
\end{proof}
\par
\medskip
Our next result shows that the polynomials $f _{6}(X)$, $f _{12}(X)$, 
$f _{18}(X)$, $f _{30}(X)$, $f _{36}(X)$, $f _{54}(X)$, $f _{84}(X)$ 
and $f _{90}(X)$ are irreducible over $\mathbb Q$. 
\par
\medskip
\begin{prop}
\label{prop3.2} 
The polynomial $f _{m}(X + 1)$ is $3$-Eisensteinian relative to the 
ring $\mathbb Z$ of integers, if $m = 3^{k} + 3^{l}$, for some 
positive integers $k$ and $l$; in particular, this holds, for $m = 6$, 
$12, 18, 30, 36, 54, 84, 90$.
\end{prop}
\par
\medskip
\begin{proof}
Note first that the free term of the polynomial $t _{m}(X) = f _{m}(X 
+ 1)$ is divisible by $3$ but is not divisible by $9$. Indeed, this 
term is equal to $2 ^{m} + 2$, and since $6 \mid m$, we 
have $2 ^{m} \equiv 1 ({\rm mod} \ 9)$ and $2 ^{m} + 2 \equiv 3 ({\rm 
mod} \ 9)$. Therefore, using Lemma \ref{lemm2.1}, one sees that it 
suffices to show that the coefficient, say $a$, of the monomial $X 
^{3^{l}}$ in the reduced presentation of $t _{m}(X)$ is divisible by 
$3$. The proof of this fact offers no difficulty because $a = (2 
^{3^{k}} + 1)\binom{m}{3 ^{l}}$ (the binomial coefficient 
$\binom{m}{3^{l}}$ is a positive integer not divisible by $3$ whereas $2 
^{3^{l}} + 1$ is divisible by $3$).
\end{proof}
\par
\medskip
Our next result gives an affirmative answer to (1) in the special 
case where $m$ is a $2$-primary number. It proves the validity of 
Theorem \ref{theo1.1}, under the condition that $m \in \{2, 4, 8, 16, 
32, 64\}$.
\par
\medskip
\begin{prop}
\label{prop3.3} 
For any $k, n \in \mathbb N$, {\rm gcd}$(f _{m}(X), f _{3n}(X)) 
= 1$, where $m = 2^{k}$.
\end{prop}
\par
\medskip
\begin{proof}
Our argument relies on the fact that $f _{m}(X) = 2g 
_{m}(X)$, $g _{m}(X)$ being a polynomial with integer coefficients, 
such that $g_{m}(X) - X^{m} - X^{m/2} - 1 = 2h _{m}(X)$, for some $h 
_{m}(X) \in \mathbb Z[X]$. This ensures that if $\rho $ is a complex 
root of $g_{m}$, $K =\mathbb Q(\rho )$ and $O_{K}$ is the ring of 
algebraic integers in $K$, then the coset $\rho + P$ is a cubic root 
of unity in the field $O/P$, for any prime ideal $P$ of $O_{K}$ of 
$2$-primary norm (i.e. a prime ideal, such that $2 \in P$). The same 
holds whenever K’/K is a finite extension, $O_{K'}$ is the ring of 
algebraic integers in $K ^{\prime }$, and $P ^{\prime }$ is a prime 
ideal in $O_{K'}$ of $2$-primary norm. The noted property of $\rho $ 
indicates that $f _{n}(\rho ) \equiv 3 \equiv 1 (mod P’)$ in case $n \in 
\mathbb N$ is divisible by $3$, which proves the non-existence of a 
common root of $f _{m}(X)$ and $f _{n}(X)$, as claimed. 
\end{proof}
\par
\medskip
The following statement provides an affirmative answer to (1), under 
the hypothesis that $m/2$ or $m/4$ is an odd primary number not 
divisible by $3$.
\par
\medskip
\begin{prop}
\label{prop3.4} 
For any prime number $p > 3$ and each pair of positive integers $k, 
n$, we have gcd$(f _{2p^{k}}(X), f _{3n}(X)) = {\rm gcd}(f 
_{4p^{k}}(X), f _{3n}(X)) = 1$.
\end{prop}
\par
\medskip
\begin{proof}
We proceed by reduction modulo $p$. Then $\bar f _{2p^{k}}(X) = \bar 
f _{2}(X) ^{p} = 2 ^{p}(X ^{2} + X + 1) ^{p}$ and $\bar f 
_{4p^{k}}(X) = \bar f _{4}(X) ^{p} = 2 ^{p}(X ^{2} + X + 1) ^{2p}$. 
This indicates that if $\hat \rho $ is a root of $\bar f 
_{2p^{k}}(X)$ or $\bar f _{4p^{k}}(X)$ in $(\mathbb Z/p\mathbb 
Z)_{{\rm sep}}$, then $\hat \rho $ is a cubic root of unity. 
Therefore, it is easily verified that $\bar f_{3n}(\hat \rho ) = -3 
\neq 0$, provided that $n$ is odd. When $n$ is even, one obtains 
similarly that $\bar f_{3n}(\hat \rho ) = 3 \neq 0$. These calculations 
prove that gcd$(\bar f_{2p^{k}}(X), \bar f_{3n}(X)) = {\rm gcd}(\bar 
f_{4p^{k}}(X), \bar f_{3n}(X)) = 1$, for each $n \in \mathbb N$. Our 
conclusion means that gcd$(\bar f_{2p^{k}}(X).\bar f_{4p^{k}}(X), 
\bar f_{3n}(X)) = 1$, which can be restated by saying that  
$u(X)f_{2p^{k}}(X)f_{4p^{k}}(X) + v(X)f_{3n}(X) = 1 + pw(X)$, for 
some $u(X), v(X), w(X) \in \mathbb Z[X]$. Suppose now that $f 
_{3n}(\beta ) = f_{2p^{k}}(\beta )f_{4p^{k}}(\beta ) = 0$, for some 
$\beta \in \mathbb C$, put $O = \{r \in \mathbb Q\colon 2^{n(r)}r 
\in \mathbb Z$, for some integer $n(r) \ge 0\}$, and denote by $O 
_{\mathbb Q(\beta )} ^{\prime }$ the integral closure in $\mathbb 
Q(\beta )$ of the ring $O$. Since $2\beta $ is an algebraic integer 
and $1 + pw(\beta ) = 0$, one obtains consecutively that $\beta \in 
O _{\mathbb Q(\beta )} ^{\prime }$ and $p$ is an invertible element 
of $O _{\mathbb Q(\beta )} ^{\prime }$; in particular, this requires 
that $1/p \in O _{\mathbb Q(\beta )} ^{\prime }$. Since, however, 
$O$ is an integrally closed subring of $\mathbb Q$, the obtained 
result leads to the conclusion that $1/p \in O$ which is not the 
case. The obtained contradiction is due to the assumption that $f 
_{3n}(X)$ and $f _{2p^{k}}(X).f_{4p^{k}}(X)$ have a common root, so 
Proposition \ref{prop3.4} is proved.
\end{proof}
\par
\medskip
Let $\bar h(Z) \in (\mathbb Z/q\mathbb Z)[Z]$ be the cubic polynomial 
defined so that $\bar h(X + X ^{-1}) = \bar g _{8}(X)/X ^{3}$, $\bar 
g _{8}(X) \in (\mathbb Z/q\mathbb Z)[X]$ being the reduction of $g 
_{8}(X) \in \mathbb Z[X]$ modulo a prime number $q > 2$ not dividing 
the discriminant $d(h)$. It is not difficult to see that the 
discriminant $d(\bar h)$ is a non-square in $(\mathbb Z/q\mathbb Z) 
^{\ast }$ if and only if $\bar h(X)$ has a unique zero lying in 
$\mathbb Z/q\mathbb Z$. When this holds, $\bar g _{8}(X)$ decomposes 
over $\mathbb Z/q\mathbb Z$ into a product of three (pairwise 
relatively prime) quadratic polynomials irreducible over $\mathbb 
Z/q\mathbb Z$. For example, this applies to the case where $q = 5$ or 
$q = 7$, which is implicitly used for simplifying the proofs of the 
following two statements.
\par
\medskip
\begin{prop}
\label{prop3.5} 
The polynomials $f _{8.5^{k}}(X)$ and $f _{n}(X)$ satisfy the equality 
{\rm gcd}$(f _{8.5^{k}}(X), f _{n}(X)) = 1$, for each $k \in \mathbb 
N$, and any $n \in \mathbb N$ divisible by $3$ and not congruent to 
$6$ modulo $24$. 
\end{prop}
\par
\medskip
\begin{proof}
It is easily verified that $2 + \sqrt{2}$, $2 - \sqrt{2}$, $1 + 
2\sqrt{2}$, $1 - 2\sqrt{2}$, $-2 - 2\sqrt{2}$ and $-2 + 2\sqrt{2}$ 
are pairwise distinct roots of $\bar f _{8}(X)$, the reduction of $f 
_{8}(X)$ modulo $5$. These roots are contained in a field $\mathbb F 
_{25}$ with $25$ elements. None of them is a primitive cubic root of 
unity: $(3 + \sqrt{2})^{3} = (2 + \sqrt{2})^{3} = -\sqrt{2}$; $(1 + 
2\sqrt{2})^{3} = 2\sqrt{2}$, $(2 + 2\sqrt{2})^{3} = 1$; $(-2 - 
2\sqrt{2})^{3} = -1$, $(-1 - 2\sqrt{2})^{3} = -2\sqrt{2}$. In other 
words, the noted elements are roots of $\bar g _{8}(X)$, the 
reduction modulo $5$ of the polynomial $g_{8}(X) = (1/2)f _{8}(X)/(X 
^{2} + X + 1)$ ($g_{8}(X) \in \mathbb Z[X]$ is irreducible over 
$\mathbb Q$). Observe that the latter two roots of $\bar g_{8}(X)$ 
are primitive $6$-th roots of $1$, whereas the remaining roots of 
$\bar g_{8}(X)$ are generators of the multiplicative group $\mathbb F 
_{25} ^{\ast }$ of $\mathbb F_{25}$. Taking further into account that 
the elements $a\sqrt{2}$, $a \in \mathbb F _{5} ^{\ast }$, are all 
primitive $8$-th roots of unity in $\mathbb F _{25}$ ($\mathbb F 
_{5}$ is the prime subfield of $\mathbb F_{25}$), and $f _{8}(X) = 
2(X^{2} + X + 1)g_{8}(X)$, one concludes that $\bar f _{8}(X)$ and 
$\bar f _{3n}(X)$ do not possess a common root, for any odd $n 
\in \mathbb N$. These calculations yield gcd$(\bar f_{8.5^{k}}(X), 
\bar f _{3n}(X)) = 1$ which allows to deduce by the method of proving 
Proposition \ref{prop3.4} that gcd$(f _{8.5^{k}}(X), f _{3n}(X)) = 1$ 
whenever $n$ is odd, and also, in the following two cases: $4 \mid 
n$; $n \equiv 6 ({\rm mod} \ 8)$. Thus Proposition \ref{prop3.5} is 
proved.
\end{proof}
\par
\medskip
\begin{prop}
\label{prop3.6} 
The polynomials $f _{8.7^{k}}(X)$ and $f _{n}(X)$ satisfy the equality
{\rm gcd}$(f _{8.7^{k}}(X), f _{n}(X)) = 1$ whenever $k$ and $n \in 
\mathbb N$, and $n$ is divisible by $6$.
\end{prop}
\par
\medskip
\begin{proof} The reductions $\bar f _{8.7^{k}}(X)$ and $\bar 
f_{8}(X)$ modulo $7$ satisfy the equality $\bar f_{8.7^{k}}(X) = 
\bar f_{8}(X)^{7^{k}}$, so it is sufficient to show $\bar f_{8}(X)$ 
and $\bar f_{n}(X)$ do not possess a common root in $(\mathbb 
Z/7\mathbb Z)_{sep}$. Our argument relies on the fact that $\sqrt{-1} 
\notin \mathbb Z/7\mathbb Z$, and $3 + \sqrt{-1}$, $3 - \sqrt{-1}$, 
$1 + 2\sqrt{-1}$, $1 - 2\sqrt{-1}$, $-2 -2\sqrt{-1}$ and $-2 + 
2\sqrt{-1}$ are all roots in $(\mathbb Z/7\mathbb Z)_{sep}$ of the 
reduction $\bar g_{8}(X)$ of $g_{8}(X)$ modulo $7$. This ensures 
that, for each of these roots, say $\rho $, $\bar f _{n}(\rho ) = 
\rho _{1} + \rho _{2} + 1$ whenever $n \in \mathbb N$ is fixed and 
divisible by $6$. Here $\rho _{1}$ and $\rho _{2}$ are $8$-th roots 
of unity in $(\mathbb Z/7\mathbb Z)_{sep}$ depending on $n$ and 
$\rho $. We show that $\bar f _{n}(\rho ) \neq 0$. Consider an 
arbitrary primitive $8$-th root of unity $\varepsilon \in (\mathbb 
Z/7\mathbb Z)_{sep}$. It is easily verified that $\varepsilon \in 
(\mathbb Z/7\mathbb Z)(\sqrt{-1}) \setminus \mathbb Z/7\mathbb Z$. 
More precisely, one obtains by straightforward calculations that 
$\varepsilon = -2(\varepsilon _{1} + \varepsilon _{2}\sqrt{-1})$, for 
some $\varepsilon _{j} \in \{-1, 1\}$, $j = 1, 2$. It is now easy to 
see that $\bar f _{n}(\rho ) \neq 0$, as claimed. Thus the assertion 
that gcd$(\bar f _{8.7^{k}}(X), \bar f _{n}(X)) = 1$, for every 
admissible pair $k, n$, becomes obvious, which completes the proof of 
Proposition \ref{prop3.6}.
\end{proof}
\par
\medskip
Proposition \ref{prop3.6} and our next result prove that gcd$(f 
_{56}(X), f_{n}(X)) = 1$, for each $n \in \mathbb N$ divisible by 
$3$. This, combined with Proposition \ref{prop3.4} and Corollaries 
\ref{coro5.6} and \ref{coro5.7}, proves the validity of Theorem 
\ref{theo1.1} in the special case where $m$ is an even biprimary 
number (see also Remark \ref{rema4.2} for the case of $m = 72$).
\par
\medskip
\begin{prop}
\label{prop3.7}
The polynomials $f_{2^{k}q}(X)$ and $f_{n}(X)$ are relatively prime, 
provided that $k \in \mathbb N$, $q \in \{5, 7\}$ and $n \in \mathbb 
N$ is odd and divisible by $3$.
\end{prop}
\par
\medskip
\begin{proof}
Denote by $\bar f_{2^{k}q}(X)$ the reduction of $f _{2^{k}q}(X)$ 
modulo $2$. It is not difficult to see that $\bar f _{2^{k}q}(X) = 
\bar f_{q}(X)^{2^{k}}$, where $\bar f_{q}$ is the reduction of 
$f_{q}(X)$ modulo $2$. This ensures that if $\alpha\in \mathbb C$ is 
an algebraic integer with $f_{2^{k}q}(\alpha ) = 0$, then $\bar f 
_{q}(\hat \alpha ) = 0$, where $\hat \alpha $ is the residue class of 
$\alpha $ modulo any prime ideal $P$ of $O_{\mathbb Q(\alpha )}$, 
such that $2 \in P$. In particular, this is the case where $\alpha $ 
is a common root of $f _{2^{k}q}(X)$ and $f_{\nu }(X)$, for some odd 
$\nu \in \mathbb N$. Observing also that $\bar f_{5}(X) = X(X + 1)(X 
^{2} + X + 1)$ and $\bar f_{7}(X) = X(X + 1)(X ^{2} + X + 1)^{2}$, 
one concludes either $\hat \alpha \in \{0, -1\}$ or $\hat \alpha $ is 
a primitive cubic root of unity. The latter possibility is clearly 
ruled out, if $3 \mid \nu $. At the same time, since $\nu $ is odd, 
$0$ and $-1$ are simple roots of $\bar f_{\nu }(X)$, so it is easy to 
see that gcd$(\bar f_{2^{k}q}(X), \bar g_{\nu }(X)) = 1$, $\bar g 
_{\nu }(X)$ being the reduction modulo $2$ of the polynomial $g_{\nu 
}(X) \in \mathbb Z[X]$ defined by the equality $f _{\nu }(X) = X(X + 
1)g _{\nu }(X)$. As $0$ and $-1$ are not roots of $f _{2^{k}q}(X)$, 
the obtained result yields consecutively gcd$(f _{2^{k}q}(X), g_{\nu 
}(X)) = gcd(f_{2^{k}q}(X), f_{\nu }(X)) = 1$, as claimed. 
\end{proof}
\par
\medskip
\begin{rema}
\label{rema3.8} Let $\bar f _{63}(X)$, $\bar f _{9}(X)$, $\bar f _{70}(X)$ 
and $\bar f _{10}(X)$ be the reductions modulo $7$ of the polynomials 
$f _{63}(X)$, $f _{9}(X)$, $f _{70}(X)$ and $f _{10}(X)$, respectively. 
It is easily verified that $f _{9}(X)$ and $f _{10}(X)$ are divisible in 
$\mathbb Z[X]$ by polynomials $g _{9}(X)$ and $g _{10}(X)$ both of 
degree $6$, which are irreducible over $\mathbb Q$. One also sees 
that $\bar g _{9}(X)$ decomposes over $\mathbb Z/7\mathbb Z$ to a 
product of two cubic  polynomials irreducible over $\mathbb Z/7\mathbb 
Z$, whereas $\bar g _{10}(X)$ is presentable as a product of three 
$(\mathbb Z/7\mathbb Z)$-irreducible quadratic polynomials. As $\bar 
f _{63}(X) = \bar f _{9}(X)^{7}$ and $\bar f _{70}(X) = \bar f 
_{10}(X)^{7}$, this yields gcd$(\bar f _{63}(X), \bar f _{70}(X)) = 
1$, which implies the existence of integral polynomials 
$u(X)$, $v(X)$ and $h(X)$, such that $u(X)f _{63}(X) + v(X)f _{70}(X) 
= 1 + 7h(X)$. We prove that gcd$(f _{63}(X), f _{70}(X)) = 1$, by 
assuming the opposite. Then $\mathbb C$ contains a common root $\beta 
$ of $f _{63}(X)$ and $f _{70}(X)$, and by (2.1) (b), $\beta $ must 
be an algebraic integer with $1 + 7h(\beta ) = 0$. This requires 
that $7$ be invertible in the ring of algebraic integers in $\mathbb 
Q(\beta )$, a contradiction proving that gcd$(f _{63}(X), f _{70}(X)) 
= 1$.
\end{rema}
\par
\medskip
It is likely that one could achieve more essential progress in the 
analysis of Question (1) (up-to its full answer), by applying  
systematically other specializations of $f _{m}(X)$ and $f _{n}(X)$ 
than those used in the proof of Proposition \ref{prop3.1}. An example 
supporting this idea is provided by the proof of the following assertion.
\par
\medskip
\begin{prop}
\label{prop3.9}
Let $n$ be a positive integer different from $6$. Then $f _{6}(X)$ 
and 
$f _{n}(X)$ have no common root except, possibly, in the 
case of $n \equiv 6$ modulo $1260$.
\end{prop}
\par
\medskip
\begin{proof}
Our starting point is the fact that $f _{6}(X)$ is irreducible over 
$\mathbb Q$ (see Proposition \ref{prop3.2}); this means that $f 
_{6}(X)$ and $f _{n}(X)$ have a common root, for a given $n \in 
\mathbb N$, if and only if $f _{6}(X)$ divides $f _{n}(X)$. Note also 
that the leading term of $f_{6}(X)$ is equal to $2$, and that $f 
_{6}(X)$ is a primitive polynomial (i.e. its coefficients are 
integers and their greatest common divisor equals $1$). These 
observations show that the complex roots of $f_{6}(X)$ are not 
algebraic integers. On the other hand, by (2.1) (b), if $n$ is odd 
and $r$ is a root of $f_{n }(X)$ and $f _{6}(X)$, then $r$ must be an 
algebraic integer. Therefore, one may assume in our further 
considerations that $n$ is even. Then it follows from Proposition 
\ref{prop3.1} that if $m \in \mathbb N$ is even and $f _{m}(X)$ 
divides $f _{n}(X)$ in $\mathbb Z[X]$, then $m - 1 \mid n - 1$ and $4 
\mid n - m$. Thus it becomes clear that $f _{6}(X) \nmid f _{n}(X)$ 
except, possibly, in the case where $n \equiv 6 \ ({\rm mod} \ 20)$. In 
order to complete the proof of Proposition \ref{prop3.9}, it remains 
to be seen that if $f _{6}(X) \mid f _{n}(X)$, then $n \equiv 6 \ {\rm 
mod} \ 126$ (by Proposition \ref{prop3.1}, $5 \mid n - 6$, so the 
divisibility of $n - 6$ by $126$ would imply $n \equiv 6$ modulo 
$126.5 = 630$; hence, by the congruence $n \equiv 6 ({\rm mod} \ 
4)$, $1260 \mid n - 6$, as claimed by Proposition \ref{prop3.9}).
\par
Observe now that $f _{6}(3) = 4^{6} + 3^{6} + 1 = 4096 + 729 + 1 = 
4826 = 2.2413 = 2.19.127$. Note also that $2^{7} = 128$ is congruent 
to $1$ modulo $127$, which implies $2^{7k'} \equiv 1 ({\rm mod} \ 
127)$, for each $k' \in \mathbb N$. Now fix an integer $k \ge 0$ and 
put $S(k) = 4^{k} + 1$. It is verified by straightforward 
calculations that $S(0) = 2$, $S(1) = 5$, $S(2) = 17$, $S(3) = 65$, 
and $S(4)$, $S(5)$ and $S(6)$ are congruent to $3$, $9$ and $33$, 
respectively, modulo $127$. It is also clear that $S(l) \equiv S(k) 
({\rm mod} \ 127)$ whenever $l$ and $k$ are non-negative integers 
with $l - k$ divisible by $7$. Thus it turns out that when $k$ runs 
across $\mathbb N$, $S(k)$ may take $7$ possible values (in fact one 
value determined by the residue class of $k$ modulo $7$).
\par
The next step towards the proof of Proposition \ref{prop3.9} is to 
show that $3$ is a primitive root of unity modulo $127$. Thereafter 
(in fact, almost simultaneously) we show that if $T(k) = 4^{k} + 
3^{k} + 1$, for any integer $k \ge 0$, then $T(k)$ is divisible by 
$127$ if and only if $k \equiv 6 ({\rm mod} \ 126)$. This particular 
fact allows us to take the final step towards our proof, as it shows 
that if $k$ is even and $f _{6}(X)$ divides $f _{k}(X)$, then 
$f _{6}(3)$ divides $f _{k}(3)$, which requires that $k \equiv 6 ({\rm 
mod} \ 126)$.
\par 
It is verified by direct calculations that $3^{63} ≡ -1 ({\rm mod} \ 
127)$ (apply the quadratic reciprocity law). Direct calculations also 
show that $3^{6}$, $3^{7}$, $3^{14}$, $3^{21}$, $3^{42}$ are 
congruent modulo $127$ to $-33$, $28$, $22$ ($127$ divides $784 – 22 
= 762$), $-19$ ($28.22 = 616$ is congruent to $-19$ modulo $127$), 
$-20$, respectively. Note also that $3^{9}$ and $3^{18}$ are 
congruent modulo $127$ to $-2$ and $4$, respectively. These 
calculations prove that $3$ is a primitive root of unity modulo 
$127$, as claimed.
\par 
The noted property of $3$ means that the residue classes modulo $127$ 
of the numbers $3^{g}\colon g = 1, \dots , 126$, form a permutation of 
numbers $1, \dots , 126$. This ensures that for any $j = 1, \dots , 7$, there 
exists a unique $s(j)$ modulo $126$, such that $S(j) + 3^{s(j)}$ is 
divisible by $127$. In order to take the final step towards our 
proof, it suffices to verify that $s(j)$ is not congruent to $j$ 
modulo $126$, for any $j \neq 6$. The verification process specifies 
this as follows: $s(0) = 9$, $s(1) = 24$, $s(2) = 101$, $s(3) = 118$, 
$s(4) = 64$, $s(5) = 65$, $s(6) = 6$. The computational part of this 
process is facilitated by the observation that $17$, $-11$, $5$ and 
$16$ are congruent modulo $127$ to $144$, $3^{5} = 243$, $132 = 
2^{2}.3.11$, and $143 = 11.13$, respectively. Proposition 
\ref{prop3.9} is proved. 
\end{proof}
\par
\medskip
\section{\bf An irreducibility criterion for integral polynomials in one 
variable}
\par
\medskip
The main result of this Section attracts interest in the question of 
whether the polynomials $f _{6n}(X)$, $n \in \mathbb N$, are 
irreducible over $\mathbb Q$. It shows that this holds in several 
special cases (which, however, is crucial for the proof of Theorem 
\ref{theo1.1}).
\par
\medskip
\begin{prop}
\label{prop4.1}
Let $f(X) \in \mathbb Z[X]$ be a polynomial of degree $n > 
0$, and let $S$ be a finite set of prime numbers not dividing the 
discriminant $d(f)$ of $f$. For each $p \in S$, denote by $n_{p}$ the 
greatest common divisor of the degrees of the irreducible polynomials 
over the field with $p$ elements, which divide the reduction of 
$f(X)$ modulo $p$.  Then every irreducible polynomial $g(X) \in 
\mathbb Z[X]$ over the field $\mathbb Q$ of rational numbers is 
divisible by the least common multiple, say $\nu $, of numbers $n
_{p}\colon p \in S$; in particular, if $\nu = n$, then $f(X)$ is 
irreducible over $\mathbb Q$.
\end{prop}
\par
\medskip
\begin{proof}
It is sufficient to observe that, by Hensel's lemma, for each $p \in 
S$, there is a degree-preserving bijection of the set of $\mathbb Q 
_{p}$-irreducible polynomials dividing $f(X)$ in $\mathbb Q _{p}[X]$ 
upon the set of $(\mathbb Z/p\mathbb Z)$-irreducible polynomials 
dividing the reduction of $f(X)$ modulo $p$ (in $(\mathbb Z/p\mathbb 
Z)[X]$). One should also note that every $\mathbb Q$-irreducible 
polynomial dividing $f(X)$ in $\mathbb Q[X]$ is presentable as a 
product of $\mathbb Q _{p}$-irreducible polynomials dividing $f(X)$ 
in $\mathbb Q _{p}[X]$.
\end{proof}
\par
\medskip
Let $f(X) \in \mathbb Z[X]$ be a $\mathbb Q$-irreducible polynomial 
of degree $n$, and let $G _{f}$ be the Galois group of $f(X)$ over 
$\mathbb Q$. It is worth mentioning that if the irreducibility of 
$f(X)$ can be deduced from Proposition \ref{prop4.1}, then $n$ 
divides the period of $G _{f}$. 
\par
\medskip
\begin{rema}
\label{rema4.2}
Using Proposition \ref{prop4.1} and a computer program for 
mathematical calculations, Junzo Watanabe proved that the 
polynomials $f _{42}(X)$, $f _{60}(X)$, $f _{66}(X)$, $f _{72}(X)$, 
$f _{78}(X)$ are irreducible over $\mathbb Q$. This result, combined 
with Proposition \ref{prop3.2}, yields gcd$(f _{m}(X), f _{n}(X)) = 
1$, for every pair $m, n \in \mathbb N$ less than $100$, such that 
$m$ is odd and $6$ divides $n$. Similarly, he proved that 
$f _{88}(X) = (X^{2} + X + 1) ^{2}g_{88}(X)$, where $g _{88}(X) \in 
\mathbb Z[X]$ has degree $84$ and is irreducible over $\mathbb Q$.
The obtained result indicates that the complex roots of $g _{88}(X)$ 
are not algebraic integers, which implies gcd$(g _{88}(X), f _{n}(X)) 
= {\rm gcd}(f _{88}(X), f _{n}(X)) = 1$, for every odd $n \in \mathbb 
N$ divisible by $3$.
\par
It would be of interest to know whether the polynomials $f 
_{6n}(X)$, $n \in \mathbb N$, are $\mathbb Q$-irreducible, and 
whether this can be obtained by applying Proposition \ref{prop4.1}.
\end{rema}
\par
\medskip
\begin{coro}
\label{coro4.3}
The polynomials $f _{88}(X)$ and $f _{n}(X) \in \mathbb Z[X]$ are 
relatively prime, for each $n \in \mathbb N$ divisible by $3$ and 
less than $100$.
\end{coro}
\par
\medskip
\begin{proof}
In view of (2.1) (b) and Remark \ref{rema4.2}, one may consider only 
the case of $6 \mid n$. Then our conclusion follows Propositions 
\ref{prop3.1}, \ref{prop3.2} and Remark \ref{rema4.2}.
\end{proof}
\par
\medskip
Statements (2.1) (d), (2.2) and Remark \ref{rema4.2}, combined with 
Propositions \ref{prop3.1}, \ref{prop3.2} and Remark \ref{rema2.2}, 
lead to the conclusion that gcd$(f _{p}(X), f _{n}(X)) = 1$, for 
every $n \in \mathbb N$ with $6 \mid n$ and $n \le 100$, and for each 
prime $p < 100$. It is worth noting that there 26 prime numbers less 
than $100$. The set of primary composite odd numbers consists of $9$, 
$25$, $27$, $49$ and $81$. 
\par
\medskip
\section{\bf Polynomials of odd orders}
\par
\medskip
Our next step towards the proof of Theorem \ref{theo1.1} aims at 
showing that
\par\noindent gcd$(f _{m}(X), f _{n}(X)) = 1$, provided that $m, n 
\le 100$ and $m$ is an odd primary number. In view of the 
observations at the end of Section 4, one may consider only the case 
where $m$ is a power of a prime $p \in \{3, 5, 7\}$. This part of our 
proof relies on (2.1) (b) and the following result.
\par
\medskip
\begin{prop}
\label{prop5.1} 
Let $p$ be a prime number and $\alpha _{n}$ a root of the polynomial  
$f _{p^{n}}(X)$, for some $n \in \mathbb N$. Suppose that $\alpha _{n}$ is 
an algebraic integer and set $\varphi _{p^{n}}(X) = p^{-1}f _{p^{n}}(X)$. 
Then $\varphi _{p}(\alpha _{n}) ^{p^{(n-1)}}$ lies in the ideal 
$pO_{\mathbb Q(\alpha_{n})}$ of the ring $O_{\mathbb Q(\alpha_{n})}$ 
of algebraic integers in $\mathbb Q(\alpha _{n})$. 
\end{prop}
\par
\medskip
Proposition \ref{prop5.1} can be deduced from Lemma \ref{lemm2.3} and 
the following lemma.
\par
\medskip
\begin{lemm}
\label{lemm5.2} 
In the setting of Lemma \ref{lemm2.3}, when $n \ge 2$, the integers 
$p^{-1}\binom{p}{j_{0}}$ and $p^{-1}\binom{p^{n}}{p^{n-1}j_{0}}$ are 
congruent modulo $p$, for each $j_{0} \in \mathbb N$, $j_{0} <  p$.
\end{lemm}
\par
\medskip
\begin{proof}
It follows from the equality $C[p^{n}, pu] = \binom{p^{n-1}}{u}$, where 
$u$ is an integer with $1 \le u < p^{n-1}$, that $\binom{p^{n}}{pu} = 
\binom{p^{n-1}}{u}.u_{p,n}$, for some element $u_{p,n}$ of the local 
ring $\mathbb Z_{(p)} = \{r/s: r, s \in \mathbb Z, \ p$ {\rm does not 
divide} $s\}$, such that $u_{p,n} - 1 \in p\mathbb Z_{(p)}$. This 
enables one to obtain step-by-step that $p^{-1}\binom{p}{j_{0}} \equiv 
p^{-1}\binom{p^{ν}}{p^{ν-1}j_{0}} ({\rm mod} \ p\mathbb Z_{(p)})$, 
$\nu = 1, \dots , n$, and so to prove Lemma \ref{lemm5.2}.
\end{proof}
\par
\medskip
The proofs of the following results rely on the explicit definitions 
of the polynomials $f _{3}(X)$, $f _{5}(X)$ and $f _{7}(X)$ (see 
(2.3)). We also need Proposition \ref{prop5.1}.
\par
\medskip
\begin{coro}
\label{coro5.3} 
We have gcd$(f _{3^{m}}(X), f _{n}(X)) = 1$ whenever $m, n \in 
\mathbb N$ and $2 \mid n$. 
\end{coro}
\par
\medskip
\begin{proof}
Let $\alpha _{m} \in \mathbb C$ be a root of both $f _{3^{m}}(X)$ 
and $f _{n}(X)$, and $P _{3}$ be a maximal ideal of $O_{\mathbb 
Q(\alpha _{m})}$, such that $3 \in P _{3}$. Then $\alpha _{m} \in 
O_{\mathbb Q(\alpha _{m})}$, so it follows from Proposition 
\ref{prop5.1} and equality (2.3) (h) that $\alpha 
_{m} ^{2} + \alpha _{m} \in P _{3}$, whence, $\alpha _{m} \in 
P _{3}$ or $\alpha _{m} + 1 \in P _{3}$. On the other hand, it is 
easy to see that $f _{n}(0)$ and $f _{n}(-1)$ are integers not  
divisible by $3$, which implies $f _{n}(\alpha _{m}) \notin P 
_{3}$. Our conclusion, however, contradicts the assumption that $f 
_{n}(\alpha _{m}) = 0$, so Corollary \ref{coro5.3} is proved.
\end{proof}
\par
\medskip
\begin{coro}
\label{coro5.4} 
The equalities gcd$(f _{5^{m}}(X), f _{n}(X)) = gcd(f _{7^{m}}(X), 
f _{n}(X)) = 1$ hold, if $m, n \in \mathbb N$ and $n$ is divisible by $6$. 
\end{coro}
\par
\medskip
\begin{proof}
It is verified by straightforward calculations that $f _{n}(0) = f 
_{n}(-1) = 2$ and $f _{n}(\varepsilon _{3}) = 3$, for each primitive 
cubic root of unity $\varepsilon _{3} \in \mathbb C$. None of these 
values is divisible by $5$ or $7$, so it is not difficult to deduce 
(in the spirit of the proof of Corollary \ref{coro5.3}) from 
Proposition \ref{prop5.1} and the definitions of $f _{5}(X)$ and 
$f _{7}(X)$ that $f _{n}(X)$ has a common root neither with $f 
_{5^{m}}(X)$ nor with $f _{7^{m}}(X)$.
\end{proof}
\par
\medskip
The following result proves the equality gcd$(f _{m}(X), f _{n}(X)) = 
1$ in the case where $m$ is odd, $m < 99$, $m$ has precisely 
two different prime divisors, and $m$ is not divisible by $9$. Here 
we note that $45$, $63$ and $99$ are all odd numbers less than $100$ 
and divisible by $9$, which have exactly two different prime divisors.
\par
\medskip
\begin{coro}
\label{coro5.5} 
Let $q \in \{3, 5, 7\}$ and $p$ be a prime number different from $2$, 
$3$ and $q$. Then gcd$(f _{qp^{\nu }}(X), f _{n}(X)) = 1$ whenever 
$\nu \in \mathbb N$, $n \in \mathbb N$ and $6$ divides $qn$.
\end{coro}
\par
\medskip
\begin{proof}
This can be obtained, proceeding by reduction modulo $p$, and arguing 
as in the proofs of Corollaries \ref{coro5.3} and \ref{coro5.4}. We 
omit the details.
\end{proof}
\par
\medskip
The next two statements prove that gcd$(f _{m}(X), f _{n}(X)) = 
1$, if $m \in \{40, 80\}$.
\par
\medskip
\begin{coro}
\label{coro5.6} 
The polynomials $f _{8.5^{k}}(X)$ and $f _{n}(X)$ satisfy the equality 
\par\noindent
{\rm gcd}$(f _{8.5^{k}}(X), f _{n}(X)) = 1$ whenever $k \in \mathbb 
N$, $n \in \mathbb N$, $3 \mid n$ and $n \le 100$.
\end{coro}
\par
\medskip
\begin{proof}
Proposition \ref{prop3.5} allows us to consider only the case of $n 
\equiv 6 ({\rm mod} \ 24)$. This amounts to assuming that $n$ equals 
$6$, $30$, $54$ or $78$. Then our calculations show that $2 + 
\sqrt{2}$, $2 - \sqrt{2}$, $1 + 2\sqrt{2}$, $1 - 2\sqrt{2}$, $-2 - 
2\sqrt{2}$ and $-2 + 2\sqrt{2}$ are roots of $\bar f _{n}(X)$, which 
implies gcd$(\bar f _{8.5^{k}}(X), \bar f _{n}(X)) \neq 
1$. Since, however, $f _{n}(X)$ is $\mathbb Q$-irreducible, for each 
admissible $n$ (see Proposition \ref{prop3.2} and Remark 
\ref{rema4.2}), it follows from Proposition \ref{prop3.1} (a) that 
gcd$(f _{8.5^{k}}(X), f _{n}(X)) = 1$, as required.
\end{proof}
\par
\medskip
\begin{coro}
\label{coro5.7}
The polynomials $f _{80}(X)$ and $f _{n}(X) \in \mathbb Z[X]$ are 
relatively prime, for every $n \in \mathbb N$ divisible by $3$ and 
less than $100$.
\end{coro}
\par
\medskip
\begin{proof}
By Proposition \ref{prop3.7}, we have gcd$(f_{80}(X), 
f_{n}(X)) = 1$ in the case where $n \in \mathbb N$ is odd and 
divisible by $3$. Note further that the same equality holds, under 
the condition that $n \in \mathbb N$, $6 \mid n$ and $n < 100$. If $n > 
80$, i.e. $n = 84, 90$ or $96$, this follows from Proposition 
\ref{prop3.2} and Remark \ref{rema2.2}. When $n \le 80$, our 
assertion can be deduced from Propositions \ref{prop3.1}, 
\ref{prop3.2} and Remark \ref{rema4.2}. 
\end{proof}
\par
\medskip
Summing-up the obtained results, one concludes that Theorem 
\ref{theo1.1} will be proved, if we show that gcd$(f _{m}(X), f 
_{n}(X)) = 1$, provided that $n \in \mathbb N$ is even, $n \le 
100$, and $m \in \{45, 63, 99\}$. We achieve this goal on a 
case-by-case basis.
\par
\medskip
\begin{coro}
\label{coro5.8} 
For each even $n \in \mathbb N$, {\rm gcd}$(f _{45}(X), f _{n}(X)) = 
1$. 
\end{coro}
\par
\medskip
\begin{proof}
Let $\bar f _{n}$ be the reduction of $f _{n}$ modulo $5$, for each $n 
\in \mathbb N$. With these notation, we have $\bar f_{45}(X) = \bar 
f_{9}(X)^{5}$ and
\par\noindent
$\bar f _{9}(X) = X(X + 1)(X - 1)^{2}(X - 2)^{2}(X - 3)^{2}$. On the 
other hand, it is easily verified that $\bar f _{n}(0) = \bar 
f _{n}(-1) = 2$ and $\bar f _{n}(j) \in \{-1, 1, 3\} \subseteq 
\mathbb Z/5\mathbb Z$, $j = 1, 2, 3$.
\end{proof}
\par
\medskip
\begin{coro}
\label{coro5.9} 
The equality {\rm gcd}$(f _{63}(X), f _{n}(X)) = 1$ holds, for any even 
number $n > 0$ at most equal to $100$.
\end{coro}
\par
\medskip
\begin{proof}
Suppose first that $6 \mid n$. Then our assertion can be proved, by 
using Remark \ref{rema2.2} and combining Proposition \ref{prop3.1} 
(a) with Proposition \ref{prop3.2} and Remark \ref{rema4.2}. It 
remains to consider the case where $3 \nmid n$. When $n = 70$, our 
conclusion follows from Remark \ref{rema3.8}, and in case $n \neq 
70$, the equality gcd$(f _{63}(X), f _{n}(X)) = 1$ is implied by 
Propositions \ref{prop3.3}, \ref{prop3.4} and \ref{prop3.7}.
\end{proof}
\par
\medskip
We are now in a position to complete the proof of Theorem \ref{theo1.1}. 
Note first that every $t \in \mathbb N$ with at least $4$ distinct 
prime divisors is greater than $100$; also, $70$ is the unique natural 
number less than $100$ and not divisible by $3$, which has 
$3$ distinct prime divisors. Therefore, the preceding assertions lead 
to the conclusion that it is sufficient to prove the equality gcd$(f 
_{99}(X), f _{n}(X)) = 1$, for every even $n \in \mathbb N$, $n \le 
100$ (the question of whether gcd$(f _{99}(X), f _{n}(X)) = 1$ 
whenever $n \in \mathbb N$ and $2 \mid n$, is open). Suppose first 
that $6 \mid n$. Then the claimed equality follows from the fact the 
roots of $f _{n}(X)$ are not algebraic integers, whereas the common 
roots of $f _{99}(X)$ and $f _{n}(X)$ (if any) must be algebraic 
integers. Henceforth, we assume that $n \le 100$ and $n$ is  not 
divisible by $3$. Applying Propositions \ref{prop3.3}, \ref{prop3.4} 
and \ref{prop3.7}, one reduces the rest of the proof of Theorem 
\ref{theo1.1} to its implementation in the special case where $n = 
70$. Denote by $\bar f _{m}$ the reduction of $f _{m}$ modulo $7$, 
for each $m \in \mathbb N$. Note that $\bar f _{70}(X) = \bar f 
_{10}(X)^{7}$ and $\bar f _{10}(X) = (X ^{2} + X + 1)^{2}\bar g 
_{10}(X)$, where $\bar g _{10}(X) \in (\mathbb Z/7\mathbb Z)[X]$ is a 
degree $6$ polynomial decomposing into a product of three pairwise 
distinct quadratic polynomials lying in $(\mathbb Z/7\mathbb Z)[X]$ 
and irreducible over $\mathbb Z/7\mathbb Z$. Suppose now that gcd$(f 
_{99}(X), f _{70}(X)) \neq 1$ and fix a common root $\beta \in 
\mathbb C$ of $f _{99}(X)$ and $f _{70}(X)$. Then $\beta $ is an 
algebraic integer, and by the preceding observations on $\bar g 
_{10}(X)$, $\hat \beta ^{48} = 1$. Here $\hat \beta $ stands for the 
residue class of $\beta $ modulo some prime ideal $P$ in the ring 
$O_{\mathbb Q(\beta )}$ of algebraic integers in $\mathbb Q(\beta )$, 
chosen so that $7 \in P$. It is clear from the equality $\hat \beta 
^{48} = 1$ that $\bar f _{99}(\hat \beta ) = (\hat \beta + 1)^{3} - 
\hat \beta ^{3} - 1 = 3(\hat \beta + 1)\hat \beta $ unless $\hat 
\beta \in \{0, -1\}$. On the other hand, the assumption that $f 
_{99}(\beta ) = f _{70}(\beta ) = 0$ requires that $\bar f _{99}(\hat 
\beta ) = 0$, which is possible only in case $\hat \beta \in \{0, 
-1\}$. The obtained contradiction is due to the hypothesis that $f 
_{99}(\beta ) = f _{70}(\beta ) = 0$. Thus it follows that gcd$(f 
_{99}(X), f _{70}(X)) = 1$, which completes the proof of the equality 
gcd$(f _{99}(X), f _{n}(X)) = 1$, for each even $n \in \mathbb N$, $n 
\le 100$. Theorem \ref{theo1.1} is proved.
\par
\medskip
\section{\bf Appendix}
\par
\medskip
After posting the first version of this preprint, Junzo Watanabe 
informed me that he had confirmed the following by a program built in 
Mathematica:
\par
\medskip
(0) $f_n(x)$  is irreducible if $n = 6m$, for all $m =1, \dots ,  100$. 
\par
\medskip
(1)
$f_n(x) / f_7(x)$   is irreducible, if $n = 6m +1$, for all $m =1, \dots ,  100$. 

\par
\medskip
(2)
$f_n(x) / f_2(x)$   is irreducible, if $n = 6m +2$, for all $m =1, \dots ,  100$. 
\par
\medskip
(3)
$f_n(x) / f_3(x)$   is irreducible, if $n = 6m +3$ for all $m =1, \dots ,  100$. 
\par
\medskip
(4)
$f_n(x) / f_4(x)$   is irreducible, if $n = 6m +4$, for all $m =1, \dots ,  100$. 

\par
\medskip
(5)
$f_n(x) / f_5(x)$   is irreducible if $n = 6m +5$ for all $m =1, \dots ,  100$. 
\par
\medskip\noindent
In view of (2.1) (c), (d) and Proposition \ref{prop3.1} (a), this 
confirmation gives an affirmative answer to (1), for pairs of 
distinct positive integers at most equal to $605$. It also answers 
in the affirmative the question posed at the end of Remark 
\ref{rema4.2}, for $n = 1, \dots , 100$. 
\par
\vskip0.2truecm
\emph{Acknowledgement.} A considerable part of the research presented 
in this note was done during my visits to Tokai University, 
Hiratsuka, Japan, in 2008/2009 and 2012. I would to thanks colleagues 
at the Department of Mathematics for their hospitality. My 
thanks are due to my host professor Junzo Watanabe also for drawing my 
attention to various aspects of Conjecture (2) and Question (1).

\end{document}